\documentclass[aps,
scriptaddress,spanish,
geometry] 
{revtex4}  

\usepackage{amsmath, amssymb}
\usepackage{makeidx}
\usepackage{amsfonts}
\usepackage{ucs}
\usepackage[utf8]{inputenc}
\usepackage[usenames,dvipsnames]{pstricks}
\usepackage{subfigure}  
\usepackage{epsfig}
\usepackage{pst-grad} 
\usepackage{pst-plot} 
\usepackage[colorlinks,hyperindex]{hyperref}
\usepackage{multirow}   
\hypersetup
{
	colorlinks,%
	citecolor=black,%
	linkcolor=black,%
	urlcolor=black,%
}

\newcommand{\bs}{\bigskip}
\newcommand{\ms}{\medskip}

\begin{document}

	\title{Riccati equations of opposite torsions from the Lie-Darboux method for spatial curves and possible applications}

 \author{Paola Lemus-Basilio}
 {\email{paola.lemus@ipicyt.edu.mx http://orcid.org/0000-0002-7364-3643}
  \affiliation{IPICyT, Instituto Potosino de Investigaci\'on Cient\'{\i}fica y Tecnol\'ogica,\\
  Camino a la presa San Jos\'e 2055, Col. Lomas 4a Secci\'on, 78216 San Luis Potos\'{\i}, S.L.P., Mexico}

 \author{Haret C. Rosu}
{\email{hcr@ipicyt.edu.mx, http://orcid.org/0000-0001-5909-1945, Corresponding author.}
  \affiliation{IPICyT, Instituto Potosino de Investigaci\'on Cient\'{\i}fica y Tecnol\'ogica,\\
  Camino a la presa San Jos\'e 2055, Col. Lomas 4a Secci\'on, 78216 San Luis Potos\'{\i}, S.L.P., Mexico}
	
	
\date{6 May 2023, revised August 2023, accepted and published September 2023}

\begin{center}
Physica Scripta 98, 105230 (2023)\\ arXiv:2305.07547
\end{center}

\begin{abstract}
A novel formulation of the Lie-Darboux method of obtaining the Riccati equations for the spatial curves in Euclidean three-dimensional space is presented. It leads to two Riccati equations that differ by the sign of torsion. The case of cylindrical helices is used as an illustrative example.
Possible applications in Physics are suggested.\\

\ms

\noindent{\it Keywords}: Riccati equation; Frenet-Serret system; Lie-Darboux method; Cylindrical helix.
\end{abstract}


\vspace{2pc}

\maketitle

\section{Introduction}

In the differential geometry of spatial curves, which is also equivalent to the kinematics of trajectories in classical mechanics, a basic result was obtained by J.F. Frenet and J.A. Serret around 1850, when they introduced the Frenet-Serret system of first-order differential equations for the moving orthogonal frame of tangent, normal, and binormal vectors
\begin{eqnarray}\label{fss}
\alpha'(s)&=&\kappa(s)\beta(s)~,\\  
\beta'(s)&=&-\kappa(s)\alpha(s)+\tau(s)\gamma(s)~,\\
\gamma'(s)&=&-\tau(s)\beta(s)~,
\end{eqnarray}
where the primes stand for the derivatives with respect to the arc length $s$ of the curve, and the coefficients $\kappa$ and $\tau$
are the curvature and torsion of the curve.
This linear system of three coupled evolution equations in the arc length of the curve is of course equivalent to a third-order differential equation in the disguised form of the following fourth-order differential equation in the tangent vector
\begin{equation}\label{fourthorder}
{\rm x}^{(iv)}-\left(2\frac{\kappa'}{\kappa}+\frac{\tau'}{\tau}\right){\rm x}'''+\left(\kappa^2+\tau^2-\frac{\kappa\kappa''-2(\kappa')^2}{\kappa^2}
+\frac{\kappa'\tau'}{\kappa\tau}\right){\rm x}''+\kappa^2\left(\frac{\kappa'}{\kappa}-\frac{\tau'}{\tau}\right){\rm x}'=0~, \quad \kappa,\tau\neq 0~.
\end{equation}

\ms

A less known result, but not less important, has been obtained towards the end of the 19th century, by S. Lie and G. Darboux, who devised a method by which the spatial curves could be described by a first-order nonlinear equation of Riccati type with coefficients expressed in terms of the curvature and torsion.

\ms

The Lie-Darboux (LD) method is mentioned in some classic textbooks in differential geometry, such as Eisenhardt's treatise \cite{eisen09} and
Struik's lectures \cite{str61}.
The goal of this paper is to introduce a more general presentation of the method, which generates two Riccati equations associated with the change of sign of the torsion.

\ms

The organization of the paper is as follows: In section 2, we briefly review the standard LD method following Struik.
In section 3, we present our generalization of the LD method. Section 4 contains the illustrative example of cylindrical helices for the explicit application of the generalized method to classical Euclidean curves. In section 5, we suggest some applications ranging from biological areas to wave/quantum optics for the obtained Riccati equations of opposite torsions and the moving frames formalism.
Finally, our conclusions are stated in Section 6.

\section{The standard LD method}\label{s2}

To have an appropriate understanding of the LD method, one should first recall the so-called spherical indicatrix mapping associated to a spatial curve in three-dimensional Euclidean space $\mathbb{E}^3$ \cite{eisen09}.
This is a one-to-one mapping between a spatial curve and a corresponding curve lying on the unit sphere centered at the origin, which is constructed by joining all the extremities of the radii parallel to the positive directions of the tangents to the given spatial curve. By this means one reduces the study of spatial curves in $\mathbb{E}^3$ to the study of curves on $\mathbb{S}^2$.

\ms

For an arbitrary spherical indicatrix curve, its tangent, normal, and binormal unit vectors satisfy the Frenet-Serret linear system (\ref{fss}) and also the algebraic equation of the unit sphere $\mathbb{S}^2$
\begin{equation}\label{e1}
\alpha^2+\beta^2+\gamma^2=1~.
\end{equation}
The idea of the LD method is to turn this additive definition of $\mathbb{S}^2$
into the factored form
\begin{equation}\label{e2}
(\alpha+i\beta)(\alpha-i\beta)=(1+\gamma)(1-\gamma)
\end{equation}
and introduce the conjugate imaginary functions $w$ and $-z^{-1}$
\begin{align}\label{w1}
&w=\frac{\alpha+i\beta}{1-\gamma}=\frac{1+\gamma}{\alpha-i\beta}~,\nonumber\\
-&\frac{1}{z}=\frac{\alpha-i\beta}{1-\gamma}=\frac{1+\gamma}{\alpha+i\beta}~,\nonumber
\end{align}
from which one can obtain the tangent, normal, and binormal in rational forms in terms of the functions $w$ and $z$
\begin{equation}\label{sphere}
\alpha=\frac{1-wz}{w-z}~, \quad \beta=i\frac{1+wz}{w-z}~, \quad \gamma=\frac{w+z}{w-z}~.
\end{equation}
Likewise in the case of the FS system, one may be interested in the evolution in the arclength variable of the functions $w$ and $-1/z$.
For the derivative of $w$, we have
\begin{equation}\label{derw}
w' =\frac{\alpha'+i\beta'}{1-\gamma}+\frac{\alpha+i\beta}{(1-\gamma)^2}\gamma'=-i\kappa w+\tau\frac{i\gamma-\beta w}{1-\gamma}~,
\end{equation}
where in the last step the derivatives from the FS system have been used. The trick in order now is to get $\alpha$ from the first definition of $w$ and substitute it in the second definition. In the resulting equation, one solves for $\beta$ and substitute it in the last equation in (\ref{derw}). The surprising result of the trick is that $\gamma$ is eliminated from the equation which becomes
\begin{equation}\label{ricstruik}
w'=-i\kappa w+\frac{i}{2}\tau w^2-\frac{i}{2}\tau~,
\end{equation}
which is a Riccati equation in which the torsion provides both the free term and the coupling to the nonlinearity.
A similar algebraic calculation shows that the function $z$ satisfies the same Riccati equation as $w$, which implies that the results for $w$ apply also to $z$.

\medskip

The second part of the LD method is to obtain the parametric equations of the curve starting
from the Riccati equation (\ref{ricstruik}). The simplest way to do it is by integrating the three components of the unit
tangent vector in (\ref{sphere}). One can see that for each of the components two particular solutions, $w$ and $z$ of (\ref{ricstruik}) are needed.
The most convenient expression for these Riccati solutions is the rational form
\begin{equation}\label{rric}
w=\frac{cf_1+f_2}{cf_3+f_4}~,\quad z=\frac{df_1+f_2}{df_3+f_4}~,
\end{equation}
where $c$ and $d$ are appropriately chosen constants that should fulfill the orthogonality relations of the FS vectors
written in terms of $w$ and $z$, and the $f_i$ are functions of the arc length $s$.
The formulas for the components of the unit tangent in this framework have been first obtained by G. Scheffers \cite{gschanw}:
\begin{eqnarray}\label{scheff}
\alpha_1&=\frac{(f_1^2-f_3^2)-(f_2^2-f_4^2)}{2(f_1f_4-f_2f_3)}~,\nonumber\\
\alpha_2&=i\frac{(f_1^2-f_3^2)+(f_2^2-f_4^2)}{2(f_1f_4-f_2f_3)}~,\\
\alpha_3&=\frac{f_3f_4-f_1f_2}{(f_1f_4-f_2f_3)}~.\nonumber
\end{eqnarray}
The parametric equations of the curve are then obtained by integrating these $\alpha_i$ along the curve
\begin{equation}\label{xyz}
x(s)=\int^s \alpha_1(\sigma) d\sigma~, \quad y(s)=\int^s \alpha_2(\sigma) d\sigma~,\quad z(s)=\int^s \alpha_3(\sigma) d\sigma~.
\end{equation}
The main issue of the LD method is to obtain the four functions $f_i(s)$, namely to be able to put the Riccati solutions in the rational form given in
(\ref{rric}).

\section{A more general approach}

Instead of (\ref{e1}), we consider now the more general algebraic equation
\begin{equation}\label{3.1}
k_1^2\alpha^2+k_2^2\beta^2+k_3 ^2\gamma^2=1~,
\end{equation}
which for $k_1=k_2=k_3=k$ is the equation of the sphere of radius $k^{-1}$ in the three-dimensional Euclidean space $\mathbb{E}^3$.
Its factored form is
\begin{align*}
(k_1\alpha+ik_2\beta)(k_1\alpha-ik_2\beta)&=(1+k_3\gamma)(1-k_3\gamma)~.
\end{align*}

The function $w$ is defined in the two ways
\begin{equation}\label{w_1}
w=\frac{k_1\alpha+ik_2\beta}{1-k_3\gamma}=\frac{1+k_3\gamma}{k_1\alpha-ik_2\beta}~,
\end{equation}
and derivating the first definition, we obtain
\begin{align*}
&w'=\frac{k_1\alpha'+ik_2\beta'}{1-k_3\gamma}+\frac{(k_1\alpha+ik_2\beta)k_3\gamma'}{(1-k_3\gamma)^2}~.
\end{align*}
Substituting the derivatives of $\alpha$, $\beta$, and $\gamma$ from the Frenet-Serret system and using (\ref{w_1}) we obtain
\begin{align*}
w'&=-i\kappa\frac{k_2\alpha+ik_1\beta}{1-k_3\gamma}+\tau\frac{ik_2\gamma-k_3\beta w}{1-k_3\gamma}~.
\end{align*}
One can see that we should take $k_1=k_2=k$ to obtain 
\begin{align}\label{ecwprim}
w'&= -i\kappa w+\tau\frac{ik\gamma-k_3\beta w}{1-k_3\gamma}~.
\end{align}
On the other hand, we can also obtain $\alpha$ from the first definition of $w$ in equation (\ref{w_1})
\begin{align*}
w(1-k_3\gamma)-ik\beta=k\alpha~,
\end{align*}
so
\begin{align*}
\alpha=\frac{w(1-k_3\gamma)-ik\beta}{k}=\frac{(1-k_3\gamma)}{k}w-i\beta~.
\end{align*}
Substituting this expression of $\alpha$ in the second definition of $w$ from equation (\ref{w_1}), we obtain
\begin{align*}
w=\frac{1+k_3\gamma}{k[\frac{(1-k_3\gamma)}{k}w-i\beta]-ik\beta}=\frac{1+k_3\gamma}{(1-k_3\gamma)w-2ik\beta}~.
\end{align*}
Thus
\begin{align*}
(1-k_3\gamma)w^2-2ik\beta w=1+k_3\gamma~.
\end{align*}
So, we obtain $\beta$ as
\begin{align*}
\beta=\frac{1+k_3\gamma-(1-k_3\gamma)w^2}{-2ikw}=\frac{i}{2}\frac{(1+k_3\gamma)}{kw}-\frac{i}{2}\frac{1-k_3\gamma}{k}w~,
\end{align*}
which upon substitution in equation (\ref{ecwprim}), and after a few easy calculus steps, gives
\begin{align*}
&w'=-i\kappa w+\tau\frac{i}{2}\frac{k_3}{k}w^2+\frac{i\tau}{2}\underbrace{\frac{2k^2\gamma-k_3(1+k_3\gamma)}{k(1-k_3\gamma)}}_Q~.
\end{align*}
Taking $k_3=k$, we obtain $Q=-1$, which leads to the Riccati equation
\begin{equation}\label{ricck}
w'=-i\kappa w+\frac{i}{2}\tau w ^2-\frac{i\tau}{2}~.
\end{equation}
On the other hand, if we take $k_3=-k$, then $Q=1$, which leads to the Riccati equation
\begin{equation}\label{ricc-k}
\tilde{w}'=-i\kappa \tilde{w}-\frac{i}{2}\tau \tilde{w}^2+\frac{i\tau}{2}~.
\end{equation}
One can see that the difference between the two Riccati equations resides in the sign of the torsion. We shall come back to this issue at the end of the next section.

\bs

\section{Example: Cylindrical helices}


\bs

This is the simplest case of curves of constant slope, i.e., curves having the ratio $\kappa/\tau={\rm const}$. For such curves, the Riccati equations
are separable and can be used to show how the LD method works.

\medskip

\underline{{\bf The case $k_3=k$}}

\medskip

The Riccati equation corresponding to this case has the form
\begin{equation}\label{4.2}
\frac{dw}{ds}=\frac{i\tau}{2}(w^2-2\xi w-1)~,
\end{equation}
where $\xi=\kappa/\tau$ is a constant which is taken as a rational number, $\xi=a/b$.

Using separation of variables and integrating, one obtains
\begin{equation}\label{4.3}
\ln\left(\frac{w-w_1}{w-w_2}\right)
 =i\frac{c}{b}\,\int^s \tau(\sigma) d\sigma + \ln K
\end{equation}
where $w_1=\xi+\sqrt{\xi^2+1}$ and $w_2=\xi-\sqrt{\xi^2+1}$ are the roots of $w^2-2\xi w-1=0$, $c=\sqrt{a^2+b^2}$, and
the arbitrary integration constant has been written as $\ln K$ for convenience.
Solving for $w$, one obtains
\begin{align}\label{gsolnk}
w&=\frac{Kw_2e^{i\frac{c}{b}\,\phi(s)}-w_1}{Ke^{i\frac{c}{b}\,\phi(s)}-1}~,\quad \phi(s)=\int^s \tau(s') ds'~,
\end{align}
from where the set of functions $f_j$, $j=1,2,3,4$, is
\begin{align*}
f_1=w_2e^{i\frac{c}{b}\,\phi(s)}~,\quad f_2=-w_1~,\quad f_3=e^{i\frac{c}{b}\,\phi(s)}~,\quad f_4=-1~.
\end{align*}
This set leads to the following expressions for the three components $\alpha_i$,
\begin{align}
\alpha_1(s)&=\frac{w_2^2-1}{w_1-w_2}\frac{e^{i\frac{c}{b}\,\phi(s)}-\Big(\frac{w_1^2-1}{w_2^2-1}\Big)e^{-i\frac{c}{b}\,\phi(s)}}{2}= i\frac{a}{bc}(a-c)\sin_{k}\left(\frac{c}{b}\,\phi(s)\right)~,\nonumber\\
\alpha_2(s)&=i\frac{w_2^2-1}{w_1-w_2}\frac{e^{i\frac{c}{b}\,\phi(s)}+\Big(\frac{w_1^2-1}{w_2^2-1}\Big)e^{-i\frac{c}{b}\,\phi(s)}}{2}
=i\frac{a}{bc}(a-c)\cos_{k}\left(\frac{c}{b}\,\phi(s)\right)~,\\
\alpha_3(s)&=-\frac{1}{ \sqrt{\xi^2+1}}=-\frac{b}{c}~,\nonumber
\end{align}
where
$\sin_{k}\left(\frac{c}{b}\,\phi(s)\right)=\frac{e^{i\frac{c}{b}\,\phi(s)}-C_{k}e^{-i\frac{c}{b}\,\phi(s)}}{2i}$, $\cos_{k}\left(\frac{c}{b}\,\phi(s)\right)=\frac{e^{i\frac{c}{b}\,\phi(s)}+C_{k}e^{-i\frac{c}{b}\,\phi(s)}}{2}$,
and $C_{k}=\frac{w_1^2-1}{w_2^2-1}=\frac{a+c}{a-c}$.

Solving for $x,y,z$, according to equation (\ref{xyz}), we find,

\begin{align}\label{e19}
&x(s)=i\frac{a}{bc}(a-c)\int^s \sin_{k}\left(\frac{c}{b}\,\phi(s')\right)\, ds'~,\nonumber\\
&y(s)=i\frac{a}{bc}(a-c)\int^s \cos_{k}\left(\frac{c}{b}\,\phi(s')\right)\, ds'~,\\
&z(s)=-\int^s\frac{1}{\sqrt{\xi^2+1}}\,ds'
=-\frac{s}{\sqrt{\xi^2+1}}=-b(s/c)~,\nonumber
\end{align}
where the integration constants have been taken zero for simplicity.

\bs

It is possible to proceed further with the calculations in an easy manner only when the torsion is a constant, say $\tau(s)=\tau$ in which case $\phi=\tau s$ and we can write
\begin{align}\label{e19}
&x(s)=i\frac{a}{bc}(a-c)\int^s \sin_{k}\left(\frac{c}{b}\,\tau s'\right)ds'=-i\frac{a}{b}(a-c)\cos_{k}\left( s/c\right)~,
\nonumber\\
&y(s)=i\frac{a}{bc}(a-c)\int^s \cos_{k}\left(\frac{c}{b}\,\tau s'\right)ds'=-i\frac{a}{b}(a-c)\sin_{k}\left( s/c\right)~,
\\
&z(s)=-\frac{s}{\sqrt{\xi^2+1}}=-b(s/c)~,\nonumber 
\end{align}
where the constant torsion has been chosen as $\tau=b/c^2$ in the last step.
Using $\sin_k^2\left(s/c\right)+\cos_{k}^2\left(s/c\right)=C_k$, one finds that
\begin{equation}\label{asquare}
x^2+y^2=a^2~,
\end{equation}
which together with the expression of $z(s)$ above shows that the helix is a cylindrical one.

\bs

\underline{{\bf The case $k_3=-k$}}

\medskip

In this case, the Riccati equation (\ref{ricc-k}) reads
\begin{align}\label{ricc-k1}
\frac{d\tilde{w}}{ds}=-\frac{i\tau}{2}(\tilde{w}^2+2\xi \tilde{w}-1)~.
\end{align}
Proceeding similarly to the previous case, one can get the $\tilde{w}$ function in the form
\begin{equation}\label{ricc-k1}
\tilde{w}=\frac{K\tilde{w}_2e^{-i\frac{c}{b}\,\phi(s)}-\tilde{w}_1}{Ke^{-i\frac{c}{b}\,\phi(s)}-1}~,
\end{equation}
where the quadratic roots are now $\tilde{w}_{1,2}=-\xi\pm \sqrt{\xi^2+1}$. Notice that $\tilde{w}_{1,2}=-w_{2,1}$ that we shall use in the expressions for the functions $f_j$, $j=1,...,4$ as inferred from equation (\ref{ricc-k1}),
\begin{align*}
f_1=-w_1e^{-i\frac{c}{b}\,\phi(s)}~,\quad f_2=w_2~,\quad f_3=e^{-i\frac{c}{b}\,\phi(s)}~,\quad f_4=-1~.
\end{align*}
When these functions are used to calculate the components $\alpha_i$ for this case, they provide the same results as (\ref{e19}). Therefore the parametric coordinates obtained in the two cases $k_3=\pm k$ coincide. This fact confirms a statement of Struik \cite{str61} that the sign of torsion cannot be determined from the parametric equations of the curve and one should use the defining differential equations to set the sign ambiguity of the torsion. The fact that our generalized approach provides two Riccati equations that differ by the sign of torsion is just the solution of this sign problem. As a final point to this section, we mention that the second-order differential equation
\begin{equation}\label{efin}
u''+i\kappa u'+\frac{\tau^2}{4}u=0~
\end{equation}
is not a defining equation in the above sense, although the defining Riccati equations (\ref{ricck}) and (\ref{ricc-k}) can be obtained by the
usual logarithmic derivative change of variable $w=-(2/i\tau)u'/u$ and $\tilde{w}=(2/i\tau)u'/u$, respectively.

\bs

\section{Possible applications in Physics}

Moving to physical applications, one should consider material equivalents of the ideal geometric concept of curve, such as  filaments, ribbons, rods, and waveguides.
Firstly, we notice that the two Riccati equations that differ by the sign of the torsion could enter the discussion of the interesting phenomenon of helix hand reversal, i.e., the spontaneous switching of a single helical structure of one handedness to its mirror image, as studied in the context of tendrils of climbing plants in \cite{gt98}, but it may occur in many other biological systems, see the references in \cite{mcmg02}.

\ms

On the other hand, for applications which are closer to quantum mechanics, we outline here some results reported by Kugler and Shtrikman \cite{ks88}.
These authors noticed the important fact that if one considers scalar products of the second and third FS equations with the binormal and normal unit vectors, respectively, one obtains
\begin{equation}\label{app1}
\gamma\cdot \frac{d\beta}{ds}=\tau~, \qquad \beta\cdot \frac{d\gamma}{ds}=-\tau~,
\end{equation}
implying that $\phi=\pm\int^{\delta s}\tau ds$ gives the amount of nonorthogonality of these vectors generated by the intrinsic torsion under a displacement of $\delta s$. In other words,
the normal unit vector develops a component along the binormal and the binormal unit vector develops an opposite equal component along the normal axis in the course of motion, which means that these vectors do not undergo parallel transport. Thus the Frenet-Serret triad do not form a locally inertial frame
implying more complicated equations of motion.

However, this nonorthogonality can be eliminated by a two-dimensional rotation of angle $\phi$ of the normal and binormal vectors
\begin{align}\label{ec3} 
\begin{pmatrix}
           {\bf U}_1 \\
           {\bf U}_2
         \end{pmatrix}=\begin{pmatrix}
           \cos\phi & -\sin\phi\\
           \sin\phi &  \cos\phi
         \end{pmatrix}\begin{pmatrix}
           \beta \\
           \gamma
         \end{pmatrix}~.
\end{align}
The arclength rates of these rotated unit vectors are given by
\begin{equation}\label{app3}
\frac{d{\bf U}_1}{ds}=-\kappa \cos\phi \, \alpha~, \qquad  \frac{d{\bf U}_2}{ds}=-\kappa \sin\phi \,\alpha~
\end{equation}
and and easy calculation shows that the $U_i$ satisfy the condition of parallel transport
\begin{equation}\label{app4}
{\bf U}_i\frac{d{\bf U}_j}{ds}=0~.
\end{equation}
Consequently, the triad $\{\alpha,{\bf U}_1,{\bf U}_2\}$ forms a locally inertial frame. After a given path along the arclength of a curve, the inertial frame will differ from the FS frame by the rotation angle $\phi$, which is the classical analog of Berry's phase that in the context of pendulums in classical mechanics is known as Hannay angle.

\ms

The case of the propagation of waves through a twisted waveguide of circular cross section made of isotropic material can be used as a simple illustration of these considerations. The phase $\phi$ depends only on the parallel transport of the transverse coordinates. The center of the circular cross section
defines a curve ${\bf x}_0(s)$ when propagation along the length of the waveguide is considered and $\alpha$ is the arclength derivative of ${\bf x}_0(s)$. In the locally inertial frame, the position of any point inside the twisted waveguide can be written as
\begin{equation}\label{app5}
{\bf x}={\bf U}_1r\cos\varphi+{\bf U}_2r\sin\varphi+{\bf x}_0(s)~.
\end{equation}
and the Helmholtz equation for a scalar field $\psi$ in the neighborhood of an arbitrary point ${\bf x}$ defined as in (\ref{app5}) takes the form
\begin{equation}\label{app6}
-\bigg[\frac{\partial^2}{\partial s^2}+\frac{1}{r^2}\frac{\partial^2}{\partial \varphi^2}+\frac{1}{r}\frac{\partial}{\partial r}r \frac{\partial}{\partial r}+\widehat{D}\bigg]\psi=E\psi~,
\end{equation}
where $\varphi=\theta +\phi$ and $\widehat{D}$ is an operator whose form for constant $\kappa$ and $\tau$ is given in \cite{ks88}, but also claimed therein to be negligible under adiabatic conditions which are usually fulfilled in experiments. Thus, if we discard $\widehat{D}$, we obtain the Helmholtz equation in the `cylindrical' coordinates $s$, $\varphi$, and $r$, which can be solved with the ansatz
\begin{equation}\label{app7}
\psi=e^{im\varphi}f_m(r,s)~,
\end{equation}
differing from the true cylindrical ansatz only in the azimuthal factor for which it is $e^{im\theta}$. Thus, under adiabatic conditions, if the standard cylindrical coordinate system is used a supplementary phase $e^{-im\phi}$ will appear in the azimuthal circular harmonics.

\ms

The authors of a recent arXiv work \cite{bds23} show that winding and linking topological invariants are functionals of intrinsic geometric quantities of torsion type characterizing a variety of topological soliton configurations. They define the vector
${\bf M}=\frac{1}{\sqrt{2}}(\beta+i\gamma)\exp\left(-i(\int\tau(x)dx+2\pi \tilde{n})\right)$, where $\tilde{n}$ in the last term of the exponent is the so-called twist number \cite{mr15}, and map it to a normalized quantum state $|\psi\rangle$ in Hilbert space. The classical Berry phase can be expressed as $i{\bf M}^*\cdot d{\bf M}/dx$, where ${\bf M}^*$ is the complex conjugate of ${\bf M}$, which corresponds through the mapping to the quantum Berry phase $i\langle\psi|\partial\psi/\partial x\rangle$.

\ms

Finally, another interesting wave phenomenon is the Gouy phase, a fundamental phase anomaly of $n\pi/2$ when a convergent light beam passes through its focus,   where $n=1$ for cylindrical beams and $n=2$ for spherical waves. The physical origin of this phase shift is still debatable despite its centenary history. There are proposals to interpret it as a Berry phase, e.g., in \cite{subb95} it is discussed for paraxial Gaussian beam optics whose propagation is described in terms of the complex beam width parameter evolving under a bilinear M\"obius transform as in (\ref{rric}). Very recently, the Gouy phase shift has been discussed in a pure quantum optical approach \cite{hbof22}, i.e., expressing wave modes in terms of photon number states. In the latter paper, the authors conclude that their results suggest a possible link between an $N$-photon state and the $N$th harmonic of a classical field, which introduces an increase of the mode order and decrease of the beam waist, in addition to doubling the frequency.

\bigskip

\section{Conclusion}

\medskip

We have revisited the Lie-Darboux method in the differential geometry of $\mathbb{E}^3$ curves. A more general approach than the one in the classic textbooks is worked out and illustrated with the example of cylindrical helices. The new feature is that two Riccati equations that differ in the sign of the torsion are obtained instead of just one as in the case of the standard method. Possible physical applications of this finding have been suggested.

\bigskip
\bigskip

\noindent {\bf Acknowledgment}\\

\noindent The authors wish to thank the reviewers for useful comments and suggestions. The first author acknowledges  CONAHCyT for the PhD fellowship.

\bigskip
\bigskip

%
%

%
%
%
%
%

\bigskip
\bigskip

\noindent {\bf ORCID iDs}\\

\noindent P. Lemus-Basilio http://orcid.org/0000-0002-7364-3643

\noindent H.C. Rosu http://orcid.org/0000-0001-5909-1945

%
%
%
%
%


\begin{thebibliography}{10}

\bibitem{eisen09} Eisenhardt L P 1909
{\em A Treatise on the Differential Geometry of Curves and Surfaces} (Boston: Ginn and Company)

\bibitem{str61} Struik D J 1961
{\em Lectures on Classical Differential Geometry} 2nd Ed. New York: Dover Publications Incorporation)

\bibitem{gschanw} Scheffers G 1923 reprinted in 2019
{\em Anwendung der Differential- und Integralrechnung auf Geometrie/Einf\"uhrung in die Theorie der Kurven in der Ebene und im Raume}
(Berlin and Leipzig: de Gruyter)   978-3-11-106496-3 (ISBN)

\bibitem{gt98} Goriely A and Tabor M 1998 Spontaneous helix hand reversal and tendril perversion in climbing plants {\em Phys. Rev. Lett.} {\bf 80} 1564
 https://doi.org/10.1103/PhysRevLett.80.1564

\bibitem{mcmg02} McMillen T and Goriely A 2002 Tendril perversion in intrinsically curved rods {\em J. Nonlin. Sci.} {\bf 12} 241
 https://doi.org/10.1007/s00332-002-0493-1

\bibitem{ks88} Kugler M and Shtrikman S 1988 Berry's phase, locally inertial frames, and classical analogous {\em Phys. Rev.} {\bf D 37} 934
 https://doi.org/10.1103/PhysRevD.37.934

\bibitem{bds23} Balakrishnan R, Dandoloff R and Saxena A 2023 Intrinsic twisted geometry underlying topological invariants arXiv: 2304.06240
 https://doi.org/10.48550/arXiv.2304.06240

\bibitem{mr15} Moffatt H K and Ricca R L 1992 Helicity and the C\u alug\u areanu invariant {\em Proc. Roy. Soc. London A} {\bf 439} 411
 https://doi.org/10.1098/rspa.1992.0159

\bibitem{subb95} Subbarao D 1995 Topological phase in Gaussian beam optics {\em Opt. Lett.} {\bf 20} 2162
 https://doi.org/10.1364/ol.20.002162

\bibitem{hbof22} Hiekkam\"aki M, Barros R F, Ornigotti M and Fickler R 2022 Observation of the quantum Gouy phase {\em Nature Photonics} {\bf 16} 828
 https://doi.org/10.1038/s41566-022-01077-w

\end{thebibliography}
\end{document}